\documentclass[11pt,leqno]{article}
\usepackage{amsfonts}
\usepackage{amsmath}
\usepackage{amssymb}
\newcommand{\pref}[1]{(\ref{#1})}
\newtheorem{theo}{Theorem}[section]
\newtheorem{lem}[theo]{Lemma}
\newtheorem{prop}[theo]{Proposition}

\title{\Large \bf {Hyper-Hermitian quaternionic K\"ahler manifolds }}
\author{{\sc Bogdan Alexandrov}
\thanks{Supported by SFB 288 "Differential geometry and quantum
physics" of DFG	and The European Contract Human Potential Programme,
Research Training Network HPRN-CT-2000-00101}
}
\date{}
\begin{document}
\maketitle
\vspace{5mm}

\begin{abstract}
We call a quaternionic K\"ahler manifold with non-zero scalar curvature, whose
quaternionic structure is trivialized by a hypercomplex structure, a
hyper-Hermitian quaternionic K\"ahler manifold. We prove that every locally
symmetric hyper-Hermitian quaternionic K\"ahler manifold is locally isometric
to the quaternionic projective space or to the quaternionic hyperbolic space.
We describe locally the hyper-Hermitian quaternionic K\"ahler manifolds with
closed Lee form and show that the only complete simply connected such
manifold is the quaternionic hyperbolic space.
\\[10mm]
\\[10mm]
\end{abstract}

\section{Introduction }

A $4n$-dimensional ($n>1$) Riemannian manifold is quaternionic K\"ahler if its
holonomy group is contained in $Sp(n)Sp(1)$. On every such manifold the bundle
of endomorphisms of the tangent bundle has a parallel $3$-dimensional
subbundle, denoted by $S^2 H$ (see Sections~\ref{sec2},\ref{sec3}), which is
locally trivialized by a triple of orthogonal almost complex structures
satisfying the quaternionic identities. Every quaternionic K\"ahler manifold is
Einstein and Alekseevsky \cite{A1} has proved that its curvature tensor has a
form which resembles that of a $4$-dimensional oriented self-dual Einstein
Riemannian manifold. This similarity allows the extension of many constructions
and results about self-dual Einstein manifolds to quaternionic K\"ahler
manifolds. Because of this, a $4$-dimensional quaternionic K\"ahler manifold is
defined to be an oriented self-dual Einstein Riemannian manifold.

There is a series of possible additional structures on a quaternionic K\"ahler
manifold: Salamon \cite{S2} has shown that $S^2 H$ always has local
sections which are complex structures; Alekseevsky, Marchiafava and Pontecorvo
\cite{P,AMP} have studied quaternionic K\"ahler manifolds with a global section
of $S^2 H$ which is almost complex or complex structure; they have proved
\cite{AM2,AMP} that if on a compact quaternionic K\"ahler manifold $S^2 H$ is
trivialized globally by an almost hypercomplex structure, then the scalar
curvature is zero, that is, the manifold is locally hyper-K\"ahler.

In the present paper we study  the quaternionic K\"ahler manifolds on which
$S^2 H$ is trivialized by a hypercomplex structure. These manifolds are
simultaneously quaternionic K\"ahler and hyper-Hermitian, so we call them {\it
hyper-Hermitian quaternionic K\"ahler (hHqK)} manifolds. To avoid the situation
of locally hyper-K\"ahler  manifolds, we require in addition that the scalar
curvature is not zero.

The simplest examples of hHqK manifolds are the the quaternionic hyperbolic
space $\mathbb{H}H^n$ and the domain of non-homogeneous quaternionic
coordinates on the quaternionic projective space  $\mathbb{H}P^n$. The whole
$\mathbb{H}P^n$ cannot be hHqK because it does not admit any almost complex
structure \cite{M}. In fact, there are no complete hHqK manifolds with positive
scalar curvature. This follows from the above mentioned result of Alekseevsky
and Marchiafava \cite{AM2} (since every complete quaternionic K\"ahler manifold
with positive scalar curvature is compact) or, alternatively, from Theorem~6.3
in \cite{S}. On the other hand, it is conjectured in \cite{AMP} that the only
complete simply connected hHqK manifold with negative scalar curvature is
$\mathbb{H}H^n$.

Further examples, generalizing the above two, are the Swann bundles
\cite{Sw}. These are principal $\mathbb{H}^* / {\mathbb{Z}_2}$-bundles over a
quaternionic K\"ahler base and have quaternionic K\"ahler metrics and a
pseudo-hyper-K\"ahler metric with hyper-K\"ahler potential, which share the
same quaternionic structure.

It is well-known that all underlying complex structures of a hyper-Hermitian
structure have the same Lee form. The condition that a  quaternionic K\"ahler
manifold is hHqK can be expressed as a differential equation for this form (see
\cite{AMP} or Proposition~\ref{prop2} below).

The hHqK structures on  $\mathbb{H}P^n$, $\mathbb{H}H^n$ and
the Swann bundles all have exact Lee forms. On the other hand, Apostolov and
Gauduchon \cite{AG} have classified locally the $4$-dimensional hHqK manifolds
with non-closed Lee form, which in addition have an orthogonal complex
structure that is not a section of
$S^2 H$: Every such manifold is locally isometric to $\mathbb{R}_+ \times S^3$
with one of the Pedersen-LeBrun metrics \cite{Ped,L}. More generally, Calderbank
\cite{C} has shown that every $4$-dimensional hHqK manifold $M$ with non-vanishing Lee form is locally of
the form $\mathbb{R} \times B$, where $B$ is a $3$-dimensional hyper-CR
Einstein-Weyl space and the metric on $M$ is explicitly given by the geometry
of $B$ (in particular, the projection $\pi : M \longrightarrow B$ is a conformal
submersion).

The two main goals of the present paper are to describe locally the hHqK
manifolds, which

{\bf A.} are locally symmetric.

{\bf B.} have closed Lee form.

With respect to problem {\bf A}, we prove in Theorem~\ref{th1} that every
locally symmetric hHqK manifold is locally homothetic to
$\mathbb{H}P^n$ or $\mathbb{H}H^n$, thus giving a positive answer to a
question of Alekseevsky and Marchiafava \cite{AM1}. The idea of the proof is
to find, in addition to the above mentioned equation for the Lee form, a
differential equation for its exterior differential and then differentiate it
until enough equations are obtained, so that the curvature tensor can
be determined. In dimension 4 Theorem~\ref{th1} is a direct consequence of a
result of Eastwood and Tod \cite{ET} about Einstein-Weyl structures on
locally-symmetric manifolds (see also \cite{AG}).

It is well-known that every $4$-dimensional hHqK manifold with closed
Lee form is locally homothetic to $S^4 \cong \mathbb{H}P^1$ or $\mathbb{R}H^4
\cong \mathbb{H}H^1$. Thus, it is enough to consider problem {\bf B} in
dimension $4n$ with $n>1$.

As already mentioned, the hHqK structures on the Swann bundles have closed
Lee forms. In Theorem~\ref{th4} we show that, conversely, every hHqK manifold,
whose Lee form is closed and has non-constant length, is locally isometric to a
Swann bundle. The proof relies on the close relation between this type of hHqK
structures and pseudo-hyper-K\"ahler metrics with hyper-K\"ahler potential,
given in Theorem~\ref{th2}.

In the remaining case of hHqK manifolds with closed Lee form of constant length
the scalar curvature is necessarily negative. Such manifolds are constructed in Theorem~\ref{th6}
in a way which resembles the construction of the Swann bundles: They are
$\mathbb{R}_+ \times \mathbb{R}^3$-bundles over a hyper-K\"ahler base.  The
converse is also true (Theorem~\ref{th7}): Every hHqK manifold with closed Lee
form of constant length is locally isometric to such a bundle.

In the last section we show that a complete simply connected hHqK manifold with
closed Lee form is homothetic to $\mathbb{H}H^n$, thus giving support to the
above mentioned conjecture of Alekseevsky, Marchiafava and Pontecorvo
\cite{AMP}.

\section{Algebraic preliminaries}\label{sec2}

Let $E$ and $H$ be the standard complex representations of $Sp(n)$ and
$Sp(1)$: $E = \mathbb{H}^n$ with $A\cdot x = Ax$ for $A \in Sp(n)$, $x \in
\mathbb{H}^n$, and $H = \mathbb{H}$ with $q\cdot y = qy$ for $q \in
Sp(1)$, $y \in \mathbb{H}$. The tensor products $E^{\otimes r} \otimes H^{\otimes s}$
are representations of $Sp(n)Sp(1) \cong Sp(n) \times_{\mathbb{Z}_2} Sp(1)$ if
$r+s$ is even, as $(-1,-1) \in Sp(n) \times Sp(1)$ acts trivially in this case.
Since $E$ and $H$ are quaternionic, these even tensor products are
complexifications of real representations of $Sp(n)Sp(1)$. For example, $E
\otimes H = T^\mathbb{C}$, where $T = \mathbb{H}^n$ with $[A,q] \in Sp(n)Sp(1)$
acting on $\xi \in \mathbb{H}^n \cong \mathbb{R}^{4n}$ by $[A,q] \cdot \xi =
A\xi \bar{q}$. This exhibits $Sp(n)Sp(1)$ as a subgroup of $SO(4n)$.

From now on, although expressing the representations of $Sp(n)Sp(1)$ in terms
of $E$ and $H$, we shall think of them as the corresponding underlying real
representations. Identifying $T$ and $T^*$ by the scalar product $g$, the
space of bilinear forms over $T$ is
\begin{equation}\label{1}
T^* \otimes T^* =
S^2 H \otimes S^2 E \oplus \mathbb{R}g \oplus \Lambda^2_0 E \oplus S^2 E \oplus
S^2 H \oplus S^2 H \otimes \Lambda^2_0 E.
\end{equation}
The first three
summands form $S^2 T^*$ and the last three form $\Lambda^2 T^*$. The space $S^2
H$ is isomorphic to the Lie algebra of $Sp(1)$. Considered as a subspace of
$End(T)$, $S^2 H = span\{I,J,K\}$, where
\begin{equation}\label{2}
I\xi = -\xi
i, \quad J\xi = -\xi j, \quad K\xi = -\xi k, \qquad \xi \in \mathbb{H}^n.
\end{equation}
The endomorphisms $I,J,K$ satisfy the quaternionic identities
\begin{equation}\label{3}
I^2 = J^2 = K^2 = - \pmb{1} = IJK,
\end{equation}
where $\pmb{1}$ is the identity operator.

Let $\pmb{L} = I \otimes I + J \otimes J + K \otimes K$, considered as an
operator on $T^* \otimes T^*$. It is $Sp(n)Sp(1)$-invariant and
$\pmb{L}^2 = 2\pmb{L} + \pmb{3}$. The eigenspace of the eigenvalue
3 is $\mathbb{R}g \oplus \Lambda^2_0 E \oplus  S^2 E$ (the remaining summands
in \pref{1} form the eigenspace of the eigenvalue $-1$). We call the bilinear
forms belonging to this eigenspace {\it $\mathbb{H}$-Hermitian} since they are
characterized by the property of being Hermitian with respect to each of $I$,
$J$, $K$. The space of skew-symmetric $\mathbb{H}$-Hermitian forms is $S^2 E$;
it is isomorphic to the Lie algebra of $Sp(n)$. The space of symmetric
$\mathbb{H}$-Hermitian forms is $\Lambda^2 E = \mathbb{R}g \oplus \Lambda^2_0
E$, with $\Lambda^2_0 E$ being the space
of symmetric trace-free $\mathbb{H}$-Hermitian bilinear forms (alternatively,
as a complex space $\Lambda^2 E = \mathbb{C}\sigma_E \oplus \Lambda^2_0 E$,
where $\sigma_E$ is the $Sp(n)$-invariant symplectic form on $E$ and
$\Lambda^2_0 E$ is the space of 2-forms whose contraction with $\sigma_E$ is
zero). The projector on the space of  $\mathbb{H}$-Hermitian bilinear forms is
obviously $\frac{1}{4}(\pmb{1} + \pmb{L})$.

An algebraic curvature tensor is called {\it hyper-K\"ahler} if it has the
algebraic properties of a curvature tensor of a hyper-K\"ahler manifold, that
is, an algebraic curvature tensor which is $\mathbb{H}$-Hermitian with respect
to the first pair of arguments (and therefore also with respect to the second
pair). The space of hyper-K\"ahler curvature tensors is \cite{S2} $S^4 E
\subset S^2 (S^2 E)$.

Let $\pi$ and $\pi_h$ be the projections on the spaces of algebraic curvature
tensors and hyper-K\"ahler algebraic curvature tensors respectively. We need
the explicit forms of $\pi$ and $\pi_h$ only in some special cases, which we
list below.

Let $R \in {T^*}^{\otimes 4}$, $\Phi \in T^* \otimes T^*$. We define $\Pi R,
\tau R, c(\Phi,R) \in {T^*}^{\otimes 4}$ by
$$\Pi R(X,Y,Z,W) = R(X,Y,Z,W) + R(Y,X,W,Z) - R(W,Z,X,Y) - R(Z,W,Y,X),$$
$$\tau R(X,Y,Z,W) = R(Y,Z,X,W),$$
$$c(\Phi,R)(X,Y,Z,W) = R(X,Y,Z,FW) = \Phi(W,R(X,Y)Z),$$
where $FX$ and $R(X,Y)Z$ are defined by
\begin{equation}\label{7} 
g(FX,Y) = \Phi (X,Y), \quad g( R(X,Y)Z,W) = R(X,Y,Z,W).
\end{equation}
If $R$ is skew-symmetric with respect to the first two arguments and satisfies
the Bianchi identity with respect to the first three arguments (that is,
$(\pmb{1} + \tau + \tau^2)R = 0$), then
$\pi R = \frac{1}{4}\Pi R$.
In particular, if $R$ satisfies these conditions, then $c(\Phi,R)$ satisfies
them as well, and therefore
$\pi c(\Phi,R) = \frac{1}{4}\Pi c(\Phi,R)$.
If, furthermore, $R$ is $\mathbb{H}$-Hermitian with respect to both first
pair and second pair of arguments and $\Phi$ is $\mathbb{H}$-Hermitian, then
\begin{equation}\label{10}
\pi_h c(\Phi,R) = \frac{1}{4}\Pi c(\Phi,R).
\end{equation}
For $\Phi, \Psi \in \Lambda^2 T^*$
\begin{equation}\label{11} 
\pi \tau \, \Phi \otimes \Psi = \frac{1}{12}(-2(\Phi \otimes \Psi +\Psi
\otimes \Phi ) + \Pi \tau \, \Phi \otimes \Psi )
\end{equation}
and for $\Phi, \Psi \in S^2 T^*$
\begin{equation}\label{12} 
\pi \tau \, \Phi \otimes \Psi = \frac{1}{4}\Pi \tau \, \Phi \otimes \Psi.
\end{equation}
Let $\iota_X$ denote the contraction by the vector $X$: for $R \in
{T^*}^{\otimes k}$ the tensor $\iota_X R \in {T^*}^{\otimes (k-1)}$ is defined
by
$$\iota_X R \, (X_1,\dots,X_{k-1}) = R(X,X_1,\dots,X_{k-1}).$$
It follows from \pref{11} and \pref{12} that for $\Phi, \Psi \in S^2 E$
\begin{eqnarray}\label{14}
\pi_h \tau \, \Phi \otimes \Psi \, (X,Y,Z,W) &=& \frac{1}{24}(-2(\Phi \otimes
\Psi + \Psi \otimes \Phi)(X,Y,Z,W) \\
& & + (\pmb{1} + \pmb{L})\iota_Y \iota_X \Pi \tau \, \Phi \otimes
\Psi \, (Z,W)). \nonumber
\end{eqnarray}

Finally, recall that the curvature tensor $R_0$ of the quaternionic projective
space $\mathbb{H}P^n$ is (in the above notations)
\begin{equation}\label{15} 
R_0(X,Y,Z,W) =  \frac{1}{2}(\pmb{1} + \pmb{L})\iota_Y \iota_X \Pi \tau \, g
\otimes g \, (Z,W) - 2\pmb{L} \, X^\flat \otimes Z^\flat \, (Y,W),
\end{equation}
where for a vector $X$ (resp. 1-form $\varphi$) $X^\flat$ (resp. $\varphi^\#$)
denotes the dual 1-form (resp. vector). Notice that $R_0$ is a
quaternionic K\"ahler curvature tensor, but not a hyper-K\"ahler curvature
tensor. Its scalar curvature is $s_0 =16n(n+2)$.

\section{HHqK manifolds: definition and general considerations}\label{sec3}

We begin this section with some well-known definitions and facts in order to
fix the notations.

An {\it almost Hermitian structure $(g,I)$} on a
manifold $M$ consists of a Riemannian metric $g$ and an orthogonal almost
complex structure $I$, that is, $I^2 = -\pmb{1}$ and $g(I\cdot ,I\cdot ) =
g$. The {\it K\"ahler form $\Omega_I$} and the {\it Lee form $\varphi$} of
$(g,I)$ (or of $I$ with respect to $g$) are defined by
$$\Omega_I (X,Y) = g(IX,Y), \quad \varphi(X) = \frac{1}{2} trace \{
I(\nabla_\bullet I)X \}.$$
In other words, $\varphi = \frac{1}{2} Id^* \Omega_I$, where $d^*$ is the
adjoint of the exterior differentiation and the action of $I$ on an
arbitrary 1-form $\psi$ is defined by the identification of $TM$ and $T^* M$ by
$g$, that is, $I\psi = -\psi \circ I$. Notice that different normalization
factors are more often used in the definition of the Lee form.

The almost complex structure $I$ is complex (or integrable) if and
only if
\begin{equation}\label{18}
(\nabla_{IX} I)IY = (\nabla_X I)Y.
\end{equation}
In this case $(g,I)$ is called a {\it Hermitian structure}. The structure
$(g,I)$ is {\it K\"ahler} if $I$ is parallel.

An {\it almost hypercomplex structure} on a manifold is defined by a triple
$(I,J,K)$ of almost complex structures, satisfying \pref{3}. If $g$ is a
Riemannian metric and $I$, $J$, $K$ are orthogonal  with respect to $g$, then
$(g,I,J,K)$ is an {\it almost hyper-Hermitian structure}. If  $(I,J,K)$ is a
{\it hypercomplex structure}, that is, if $I$, $J$, $K$ are integrable, the
structure $(g,I,J,K)$ is called {\it hyper-Hermitian}. It is well-known that in
this case $I$, $J$, $K$ share the same Lee form $\varphi$, which is also the Lee form of
each of the complex structures in the $S^2$-family, determined by them.

An almost hyper-Hermitian structure is {\it hyper-K\"ahler} if each of  $I$,
$J$, $K$ is parallel. Every hyper-K\"ahler manifold is Ricci flat.

An {\it almost quaternionic Hermitian structure} on a manifold $M$ consists of
a Riemannian metric and a $3$-dimensional subbundle of the bundle of
endomorphisms of $TM$, which is locally trivialized by a triple of orthogonal
almost complex structures, satisfying \pref{3}. If $\dim M = 4n$ with $n>1$,
the structure is called {\it quaternionic K\"ahler} when the subbundle is
parallel. Equivalently, a $4n$-dimensional ($n>1$) Riemannian manifold is
quaternionic K\"ahler if its holonomy group is contained in $Sp(n)Sp(1)$.

We use the same notation for a $Sp(n)Sp(1)$-representation and the
corresponding bundle, associated to the principal $Sp(n)Sp(1)$-bundle given by
the holonomy reduction. For example, the defining $3$-dimensional subbundle is
$S^2 H$.

The condition that $S^2 H$ is parallel can be expressed in terms of a local
trivializing almost hypercomplex structure $(I,J,K)$ by the equation
\begin{equation}\label{19}
\bigl( \nabla_X I, \nabla_X J, \nabla_X K \bigr) = \bigl( I,J,K \bigr)
D(a,b,c)(X), \end{equation}
where $D(a,b,c)$ is an $\mathfrak{so}(3)$-valued 1-form, given by
$$D(a,b,c) =
\begin{pmatrix}
0 & -c & b \\
c & 0 & -a \\
-b & a & 0
\end{pmatrix},
\qquad a,b,c \in \Gamma (T^* M).$$
\begin{lem}\label{lem1}
i) A $4n$-dimensional ($n>1$) almost quaternionic Hermitian manifold is
quaternionic K\"ahler if and only if the K\"ahler forms of a local
trivializing almost hypercomplex structure $(I,J,K)$ satisfy
\begin{equation}\label{21} 
\bigl( d\Omega_I, d\Omega_J, d\Omega_K \bigr) = \bigl(
\Omega_I,\Omega_J,\Omega_K \bigr) \wedge D(a,b,c).
\end{equation}
The forms $a$, $b$, $c$  coincide with those in \pref{19}.

ii) An almost hyper-Hermitian manifold is hyper-K\"ahler if and only if
$d\Omega_I = d\Omega_J = d\Omega_K = 0$.
\end{lem}

\noindent {\it Proof:} i) It follows from \pref{19} that \pref{21} is satisfied
on a quaternionic K\"ahler manifold. Conversely, if \pref{21} is satisfied,
then the algebraic ideal of the exterior algebra $\Lambda T^* M$,
generated by $S^2 H$, is a differential ideal and the fundamental form
$\Omega_I \wedge \Omega_I + \Omega_J \wedge \Omega_J + \Omega_K \wedge
\Omega_K$ is parallel. Thus, by Theorem~2.2 in \cite{Sw}, the
manifold is quaternionic K\"ahler.

That the forms $a$, $b$, $c$  coincide with those in \pref{19} follows from the
injectivity of the map
$$T^* M \oplus T^* M \oplus T^* M \ni (\alpha,\beta,\gamma) \mapsto \alpha
\wedge \Omega_I + \beta \wedge \Omega_J + \gamma \wedge \Omega_K \in \Lambda^3
T^* M$$
under the given dimension assumption.

Part ii) is proved by Hitchin \cite{H}. \hfill $\square$

\vspace{3mm}
\noindent {\bf Remark 1} Lemma~\ref{lem1} is in fact of purely algebraic
nature and therefore it is easily seen that its first part is a direct
consequence of its second part.

\vspace{3mm}
Every quaternionic K\"ahler manifold is Einstein. Following \cite{AM1,AM2}, we
denote by $\nu$ its {\it reduced scalar curvature}, $\nu = \frac{4s}{s_0} =
\frac{s}{4n(n+2)}$, where $s$ is the (constant) scalar curvature. The curvature
tensor of a quaternionic K\"ahler manifold has the form \cite{A1,S}
\begin{equation}\label{22} 
R=\frac{1}{4}\nu R_0 + R',
\end{equation}
where $R_0$ is the (parallel) curvature tensor of $\mathbb{H}P^n$, given by
\pref{15}, and $R' \in \Gamma(S^4 E)$ (that is, $R'$ is a hyper-K\"ahler
curvature tensor).

When the dimension is 4, the definition of a quaternionic K\"ahler manifold
gives nothing more than an oriented Riemannian manifold. But its curvature
tensor has the form \pref{22} only if it is self-dual and Einstein. Because of
this, we define, as is usually done, a $4$-dimensional quaternionic K\"ahler
manifold to be an oriented self-dual Einstein manifold.

The following fact is well-known.
\begin{lem}\label{lem2}
On a quaternionic K\"ahler manifold the K\"ahler forms of a local
trivializing almost hypercomplex structure $(I,J,K)$ satisfy
$$\bigl( da + b \wedge c, db + c \wedge a, dc + a \wedge b \bigr) = - \nu
\bigl( \Omega_I,\Omega_J,\Omega_K \bigr),$$
where $a$, $b$, $c$ are the 1-forms in \pref{19}.
\end{lem}

A quaternionic K\"ahler manifold with vanishing scalar curvature is locally
hyper-K\"ahler. Since we would like to avoid this situation, we assume in the
sequel that the quaternionic K\"ahler manifolds satisfy the additional
requirement of having non-zero scalar curvature.

Let $(I,J,K)$ be a local almost hypercomplex structure on a
quaternionic K\"ahler manifold, trivializing $S^2 H$. Then it follows from
\pref{19} that the Lee forms of $I$, $J$, $K$ are $-\frac{1}{2}(Jb + Kc)$,
$-\frac{1}{2}(Kc + Ia)$, $-\frac{1}{2}(Ia + Jb)$ respectively. Using also
\pref{18}, we see that $I$, $J$, $K$ are integrable if and only if $Ia
= Jb = Kc$. Thus we obtain (see also \cite{AMP,AMP2})
\begin{prop}\label{prop1}
On a quaternionic K\"ahler manifold a local
trivializing almost hypercomplex structure $(I,J,K)$ is hypercomplex if and
only if there exists a 1-form $\varphi$ such that in \pref{19} $a = I\varphi$,
$b = J\varphi$, $c = K\varphi$. In this case $\varphi$ is the common Lee form
of the hypercomplex structure.
\end{prop}

The next proposition gives the necessary and sufficient condition under which
$S^2 H$ is locally trivialized by a hypercomplex structure (see also
\cite{AMP}).
\begin{prop}\label{prop2}
The bundle $S^2 H$ on a quaternionic K\"ahler manifold is locally trivialized
by a hypercomplex structure with Lee form $\varphi$ if and only if there exists
a 2-form $\Phi \in \Gamma(S^2 E)$ such
that
\begin{equation}\label{25} 
\nabla \varphi = \frac{1}{2} \bigl( (-\pmb{1} + \pmb{L}) \, \varphi
\otimes \varphi - \nu g \bigr) + \Phi.
\end{equation}
In this case necessarily $\Phi = \frac{1}{2} d\varphi $.
\end{prop}
Notice that the operator $\pmb{L}$ on a quaternionic K\"ahler manifold is
independent of the choice of a local trivializing almost hypercomplex structure
and is parallel by \pref{19}.

\vspace{2mm}
\noindent{\it Proof:} Let $(I,J,K)$ be a local trivializing
almost hypercomplex structure. We want to find a local trivializing
hypercomplex structure $(I',J',K')$. It follows from \pref{3} that
$\bigl( I',J',K' \bigr) = \bigl( I,J,K \bigr) S$
for some $S \in SO(3)$. By \pref{19} and Proposition~\ref{prop1}, we get
\begin{equation}\label{27} 
dS + DS = SD',
\end{equation}
where $D=D(a,b,c)$, $D' = D(I'\varphi, J'\varphi,K'\varphi)$. Let
$\widetilde{D} = D(I\varphi, J\varphi,K\varphi)$.
Then $SD' = \widetilde{D} S$ and \pref{27} becomes
$$dS = (\widetilde{D} - D)S.$$
This equation has a solution locally if and only if
$$d(\widetilde{D} - D) - (\widetilde{D} - D) \wedge (\widetilde{D} - D) = 0.$$
Using Lemma~\ref{lem2}, this is equivalent to
$$\nabla_X \varphi (Y) + \nabla_{IY} \varphi (IX) = J\varphi (X) J\varphi (Y) +
K\varphi (X) K\varphi (Y) - \nu g(X,Y)$$
and the two similar equations obtained by cyclic permutations of $I$, $J$, $K$.
Then it is easy to see that these equations are equivalent to the requirement
that the symmetric part of $\nabla \varphi$ is $$\frac{1}{2} \bigl(
(-\pmb{1} + \pmb{L}) \, \varphi \otimes \varphi - \nu g \bigr)$$
and its skew-symmetric part is $\mathbb{H}$-Hermitian. \hfill $\square$

\vspace{3mm}
\noindent {\bf Definition } A {\it hyper-Hermitian quaternionic K\"ahler
(hHqK)} manifold is a quaternionic K\"ahler manifold such that $S^2 H$ is
trivialized by a hypercomplex structure. The Lee form of the hypercomplex
structure is called the {\it Lee form} of the hHqK manifold.
\begin{prop}\label{prop3}
Let $M$ be a hHqK manifold with Lee form $\varphi$, $\xi = \varphi^\#$ and
$\Phi = \frac{1}{2} d\varphi $. Then
\begin{equation}\label{33}
\nabla_Z \Phi = \iota_\xi \iota_Z R' - \varphi(Z)\Phi - \frac{1}{2}
(\pmb{1} + \pmb{L}) \, \varphi \wedge \iota_Z \, \Phi;
\end{equation}
\begin{equation}\label{34} 
\nabla_\xi R' = \nu R' + A_1 + A_2 + A_3,
\end{equation}
where
$$A_1 = - \frac{1}{2}\Pi c \bigl( (\pmb{1} + \pmb{L}) \, \varphi \otimes
\varphi, R'\bigr),$$
$$A_2 = - \Pi c(\Phi, R'),$$
$$A_3 = 12\pi_h \tau \, \Phi \otimes \Phi$$
(the right-hand side is given by \pref{14});
\begin{equation}\label{38} 
\nabla^2_{\xi,U} R' = \frac{3}{2}\nu \nabla_U R' + B_{1U} + B_{2U} + B_{3U} +
B_{4U} +  B_{5U},
\end{equation}
where
\begin{eqnarray}
{B_1}_U (X,Y,Z,W) &=&  \frac{1}{2}(-\pmb{1} + \pmb{L})\bigl(
\varphi \otimes \nabla_\bullet R'(X,Y,Z,W)\bigr)(U,\xi) \nonumber \\
& & - \frac{1}{2}\Pi c \bigl( (\pmb{1} + \pmb{L}) \, \varphi \otimes \varphi,
\nabla_U R' \bigr)(X,Y,Z,W), \nonumber
\end{eqnarray}
$${B_2}_U =  -\nabla_{FU} R' - \Pi c(\Phi, \nabla_U R')$$
($F$ is the endomorphism corresponding to $\Phi$ by \pref{7}),
$${B_3}_U = \nu \varphi (U) R' + \frac{1}{2}\nu \Pi c \bigl( (\pmb{1} +
\pmb{L}) \, \varphi \otimes U^\flat, R' \bigr),$$
$${B_4}_U = -\Pi c \bigl( (\pmb{1} + \pmb{L}) \, \iota_U \Phi \otimes \varphi,
R' \bigr) + 24\pi_h \tau \, \iota_\xi \iota_U R' \otimes \Phi,$$
$${B_5}_U = -\varphi (U) A_3 -\frac{1}{2} \Pi c \bigl( (\pmb{1} + \pmb{L}) \,
\varphi \otimes U^\flat, A_3 \bigr).$$
\end{prop}

\noindent {\it Proof:} The first equality is proved by using \pref{25} to
calculate $\nabla ^2_{X,Y} \varphi (Z)$, then antisymmetrizing with
respect to $X$ and $Y$ to get $R(X,Y,Z,\xi )$ and using \pref{22}, \pref{15} and
$d\Phi = 0$.

The second equality is proved in a similar fashion by using \pref{33} and
\pref{25} to calculate  $\nabla ^2_{Z,W} \Phi $ and then antisymmetrizing
to get $R(Z,W)\Phi$.

The equality \pref{38} is obtained by differentiating \pref{34} with respect to
$U$ and using \pref{25} and \pref{33} to substitute $\nabla_U \varphi$ and
$\nabla_U \Phi$ and \pref{34} to express $A_2$ through $\nabla_\xi R'$, $A_1$
and $A_3$.  \hfill $\square$

\vspace{3mm}
\noindent {\bf Remarks}

{\bf 2)} By \pref{25}, we can determine all components of $\nabla \varphi$
with respect to the decomposition \pref{1}. For example,
\begin{equation}\label{44}
d^* \varphi = 2n\nu - \vert \varphi \vert^2.
\end{equation}

{\bf 3)} It is clear from the proof of Proposition~\ref{prop3} that it is
derived only from  \pref{25}. This means that it remains true on the whole set
where the solution $\varphi$ of  \pref{25} is defined, although the
hyper-Hermitian structure, corresponding to  $\varphi$, may exist only on a
smaller set.

{\bf 4)} It follows from \pref{33} that $\Phi$ is co-closed and therefore
harmonic, a result obtained in \cite{AMP}.

{\bf 5)} By \pref{10} and \pref{14},
$A_1$, $A_2$, $A_3$, ${B_1}_U$,\dots ,${B_5}_U$ are hyper-K\"ahler curvature
tensors. It is also easily seen that $B_1$,\dots ,$B_5$ satisfy the second
Bianchi identity and therefore they have the algebraic properties of a
covariant derivative of a hyper-K\"ahler curvature tensor.

\section{Locally symmetric hHqK manifolds}

In this section we give a positive answer to a question of Alekseevsky and
Marchiafava \cite{AM1} concerning the symmetric quaternionic K\"ahler
manifolds, which are locally hypercomplex. In dimension 4 our theorem is a
direct consequence of a result  of Eastwood and Tod \cite{ET}.
\begin{theo}\label{th1}
A locally symmetric hHqK manifold is locally homothetic to $\mathbb{H}P^n$ or
$\mathbb{H}H^n$ and its Lee form is closed.
\end{theo}

\noindent {\it Proof:} For a vector $X$ we denote $span \, \{ X,IX,JX,KX \}$ by
$span_\mathbb{H} \{ X \}$ and the orthogonal complement of $span_\mathbb{H} \{
X \}$ by $span_\mathbb{H} \{ X \}^\bot$.

The vanishing of $\nabla R'$ means that \pref{34} and \pref{38} are
reduced to
\begin{equation}\label{45} 
0 = \nu R' + A_1 + A_2 + A_3,
\end{equation}
\begin{equation}\label{46}
0 = B_3 + B_4 +  B_5.
\end{equation}

The Lee form $\varphi$ cannot be zero on an open set (otherwise
Proposition~\ref{prop1} and Lemma~\ref{lem2} imply $\nu = 0$). Thus, it is
enough to prove that $R'$ and $\Phi$ vanish at the points where $\varphi \not =
0$. We do this in five consecutive steps. At every step we put in
\pref{46} arguments $U$, $X$, $Y$, $Z$, $W$ of certain type and prove that
certain components of $R'$ and $\Phi$ vanish.

{\it Step 1.} $U,X,Y,Z,W \in span_\mathbb{H} \{ \xi \}$.
\\
In this case it follows that $\Phi$ and $R'$ vanish if all their arguments are
in $span_\mathbb{H} \{ \xi \}$. This completes the proof if the dimension is 4.

{\it Step 2.} $U \in span_\mathbb{H} \{ \xi \}^\bot$, $X,Y,Z,W
\in span_\mathbb{H} \{ \xi \}$.
\\
Then we see that $R'$ vanishes if three of its arguments lie in
$span_\mathbb{H} \{ \xi \}$.

{\it Step 3.} $U = X = Z = \xi $, $Y,W \in
span_\mathbb{H} \{ \xi \}^\bot$.
\\
Then we get
\begin{equation}\label{401}
A_3 (\xi, Y,\xi,W) = \nu R'(\xi, Y,\xi,W).
\end{equation}
Now we put $Y = W = F\xi$ in this equation ($F\xi$ is orthogonal to
$span_\mathbb{H} \{ \xi \}$ since $F$ is an antisymmetric endomorphism commuting
with $I$, $J$, $K$). This yields
\begin{equation}\label{402}
-3\vert F\xi \vert^4 = \nu R'(\xi,F\xi,\xi,F\xi).
\end{equation}
It follows from \pref{25} and \pref{33} that
\begin{equation}\label{403}
\xi (\vert F\xi \vert^4) = -2(4\vert \varphi \vert^2 + \nu ) \vert F\xi
\vert^4,
\end{equation}
$$\xi \bigl( R'(\xi,F\xi,\xi,F\xi) \bigr) = -(5\vert \varphi \vert^2 + 2\nu )
R'(\xi,F\xi,\xi,F\xi) + 2R'(\xi,F^2 \xi,\xi,F\xi).$$
Putting in \pref{45} arguments $\xi$, $F\xi$, $\xi$, $F\xi$, we see that
$$2R'(\xi,F^2 \xi,\xi,F\xi) = (-\vert \varphi \vert^2 + 2\nu
)R'(\xi,F\xi,\xi,F\xi)$$
and therefore
\begin{equation}\label{404}
\xi \bigl( R'(\xi,F\xi,\xi,F\xi) \bigr) = -6\vert \varphi \vert^2
R'(\xi,F\xi,\xi,F\xi).
\end{equation}
Now it follows from \pref{402}, \pref{403}, \pref{404} that
$$(\vert \varphi \vert^2 + \nu ) \vert F\xi \vert^4 = 0.$$
Thus $F\xi = 0$ or $\vert \varphi \vert^2 = -\nu$ is constant. In the latter case, using \pref{25}, we obtain
$$0 =  (F\xi ) (\vert \varphi \vert^2) = -2\vert F\xi \vert^2.$$
Hence $F\xi = 0$, which means that $\Phi$ vanishes if one of its arguments
belongs to  $span_\mathbb{H} \{ \xi \}$. This implies that the same is true for
$A_3$ and by \pref{401} we get $R'(\xi, Y,\xi,W) = 0$. Hence, $R'$ vanishes if
two of its arguments are in $span_\mathbb{H} \{ \xi \}$.

{\it Step 4.} $U,X,Z \in span_\mathbb{H} \{ \xi \}^\bot$, $Y = W = \xi$.
\\
In this case it follows that $R'$ vanishes if one of its arguments lies in
$span_\mathbb{H} \{ \xi \}$.

{\it Step 5.} $U = \xi$, $X,Y,Z,W \in span_\mathbb{H} \{ \xi \}^\bot$.
\\
Then we obtain
\begin{equation}\label{52} 
A_3(X,Y,Z,W) = \nu R'(X,Y,Z,W).
\end{equation}
This is also true for arbitrary $X$, $Y$, $Z$, $W$, since if any of them
belongs to $span_\mathbb{H} \{ \xi \}$, then both $A_3$ and $R'$ vanish. Hence,
$$\nabla_\xi A_3 = \nu \nabla_\xi R' = 0.$$
But from \pref{33} we have
$\nabla_\xi \Phi = -\vert \varphi \vert^2 \Phi$,
and therefore
$\nabla_\xi A_3 = -2\vert \varphi \vert^2 A_3$.
Thus $A_3 = 0$ and by \pref{52} $R' = 0$. The vanishing of $A_3$ implies also
$\Phi = 0$.  \hfill $\square$

\section{HHqK manifolds with closed Lee form}

In this section $M$ is a hHqK manifold with closed Lee form $\varphi$ and the
hypercomplex structure is $(I,J,K)$.

Because of $\Phi = \frac{1}{2} d\varphi = 0$, \pref{25} and \pref{33} become
\begin{equation}\label{56} 
\nabla \varphi = \frac{1}{2} \bigl( (-\pmb{1} + \pmb{L}) \, \varphi
\otimes \varphi - \nu g \bigr),
\end{equation}
\begin{equation}\label{57}
R'(X,Y,Z,\xi) = 0.
\end{equation}
If the dimension is 4, \pref{57} implies $R' = 0$. Thus a $4$-dimensional
hHqK manifold with closed Lee form is locally homothetic to $\mathbb{H}P^1
\cong S^4$ or $\mathbb{H}H^1 \cong \mathbb{R}H^4$, a fact, which is well-known.
Hence, for the rest of this section we can assume that $\dim M = 4n$ with
$n>1$.

By \pref{56},
\begin{equation}\label{106}
\nabla_{\xi } \xi = -\frac{1}{2} (\vert \varphi \vert^2 + \nu ) \xi.
\end{equation}
Thus, after a change of the parameter,  the integral curves of $\xi$ are
geodesics.

From \pref{56} and Proposition~\ref{prop1}, we get
\begin{equation}\label{58} 
\nabla (I\varphi ) = \frac{1}{2} \bigl( -\varphi \otimes I\varphi - I\varphi
\otimes \varphi - J\varphi \wedge K\varphi - \nu \Omega_I \bigr)
\end{equation}
and similarly by cyclic permutations of $I$, $J$, $K$. Equation \pref{56}
also implies
\begin{equation}\label{60} 
d(\vert \varphi \vert^2 + \nu ) = - (\vert \varphi \vert^2 + \nu )\varphi.
\end{equation}
Since our considerations will be local, we can assume that $\varphi = df$ for
some function $f$. Thus from \pref{60} it follows that
\begin{equation}\label{61} 
(\vert \varphi \vert^2 + \nu )e^f = C,
\end{equation}
where $C$ is a constant.

Let $\psi = d e^f$ and $\eta = \psi ^\# = e^f \xi$. By
\pref{56} and \pref{58}, we obtain
\begin{equation}\label{62}
\nabla \psi = \frac{1}{2} e^f \bigl( (\pmb{1} + \pmb{L}) \, \varphi
\otimes \varphi - \nu g \bigr),
\end{equation}
\begin{equation}\label{63}
\nabla (I\psi ) = \frac{1}{2} e^f \bigl( \varphi \wedge I\varphi - J\varphi
\wedge K\varphi - \nu \Omega_I \bigr)
\end{equation}
and similarly by cyclic permutations of $I$, $J$, $K$.

Equation \pref{62} shows that $\nabla \psi$ is $\mathbb{H}$-Hermitian and
therefore $\eta $ is an infinitesimal quaternionic automorphism. Even more,
using also Proposition~\ref{prop1}, we see that $\eta$ is an infinitesimal
automorphism of $I$, $J$, $K$, that is, an infinitesimal hypercomplex
automorphism. It follows again from \pref{62} that, similarly to $\xi$, after a
change of the parameter the integral curves of $\eta$ are geodesics.

By \pref{63} $I\eta$, $J\eta$, $K\eta$ are Killing vector fields, which are
also infinitesimal quaternionic automorphisms (in fact, every Killing vector
field on a quaternionic K\"ahler manifold is an infinitesimal quaternionic
automorphism, see \cite{S2}). They are infinitesimal hypercomplex
automorphisms only if $\vert \varphi \vert^2 + \nu = 0$.

It follows from \pref{62}, \pref{63}, \pref{22}, \pref{15} and \pref{57} that
$span_\mathbb{H} \{ \eta \}$ is a totally geodesic distribution with integral
manifolds of constant curvature $\nu$ (larger totally geodesic quaternionic
distributions containing $span_\mathbb{H} \{ \eta \}$ exist on hHqK manifolds with closed
Lee form, see \cite{AM3,AMP}). The commutators of $\eta$, $I\eta$, $J\eta$,
$K\eta$ are given by
\begin{equation}\label{64} 
[\eta,I\eta] = 0, \quad [I\eta,J\eta] = CK\eta
\end{equation}
and the same with cyclic permutations of $I$, $J$, $K$.

Hence, if $C \not = 0$ (that is, if $\vert \varphi \vert^2 + \nu \not = 0$),
$I\eta$, $J\eta$, $K\eta$ induce an infinitesimal isometric action of $Sp(1)$
and together with $\eta$ they give rise to an infinitesimal quaternionic action
of $\mathbb{H}^*$ on $M$. This situation very much resembles the one in the
case of hyper-K\"ahler manifold with hyper-K\"ahler potential, described by
Swann \cite{Sw}. Below we show that these two situations are closely related.

We recall the definition of a hyper-K\"ahler potential in the pseudo-Riemannian
settings.

A function $\mu$ on a pseudo-K\"ahler manifold $(M,g_0,I)$ is called a {\it
K\"ahler potential} if
\begin{equation}\label{65}
\frac{1}{2}dId\mu = \Omega_I^0.
\end{equation}
A function $\mu$ on a pseudo-hyper-K\"ahler manifold $(M,g_0,I,J,K)$ is called
a {\it hyper-K\"ahler potential} if it is a K\"ahler potential for each of the
underlying pseudo-K\"ahler structures. As shown by Swann \cite{Sw}, this is
equivalent to
\begin{equation}\label{66} 
\nabla d\mu = g_0.
\end{equation}
Hence,
\begin{equation}\label{67} 
d(g_0(d\mu,d\mu)) = 2 d\mu .
\end{equation}
\begin{theo}\label{th2}
i) Let $(M,g)$ be a hHqK manifold with closed Lee form $\varphi$ and reduced
scalar curvature $\nu$, such that $g(\varphi,\varphi) + \nu \not = 0$.
Then with respect to the same hypercomplex structure
$$g_0 = \frac{1}{\nu ^2} ( g(\varphi,\varphi ) + \nu ) \left(
 (\pmb{1} + \pmb{L}) \, \varphi \otimes \varphi  + \nu g \right)$$
is a pseudo-hyper-K\"ahler metric with hyper-K\"ahler potential
$\mu = \frac{2}{\nu^2} (g(\varphi,\varphi ) + \nu)$.
The signature of $g_0$ is Riemannian when $\nu (g(\varphi ,\varphi) + \nu) > 0$
(and therefore always when $\nu >0$), and $(4,4(n-1))$ with positive sign on
$span_\mathbb{H} \{ \varphi \}$ when $\nu (g(\varphi ,\varphi) + \nu) < 0$.

ii) Let $(M,g_0)$ be a pseudo-hyper-K\"ahler manifold with hyper-K\"ahler
potential $\mu$. Then for each $p \not = 0$
$$g_p = - \frac{p}{(pg_0(d\mu,d\mu ) + 1)^2 } (\pmb{1} + \pmb{L}) \, d\mu
\otimes d\mu + \frac{1}{pg_0(d\mu,d\mu ) + 1}g_0 $$
(defined on the submanifold $pg_0(d\mu,d\mu ) + 1 \not = 0$) forms together
with the given hypercomplex structure a pseudo-hHqK structure with Lee form
$\varphi_p = - d\ln \vert pg_0(d\mu,d\mu ) + 1 \vert$ and reduced scalar
curvature $4p$. The metric $g_p$ is positive definite when $pg_0(d\mu,d\mu ) + 1
>0$ and $g_0$ is positive definite, and when $pg_0(d\mu,d\mu ) + 1 <0$ and $g_0$
has signature $(4,4(n-1))$ with positive sign on $span_\mathbb{H} \{ d\mu \}$.

iii) If $g_0$ is constructed from $g$ as in i) and $g_p$ from $g_0$ as in ii),
then $g = g_{\frac{1}{4}\nu}$. Similarly, if we start with a
pseudo-hyper-K\"ahler metric and construct $g_p$, then the
pseudo-hyper-K\"ahler metric constructed from $g_p$ is just the initial one.
\end{theo}

\noindent{\it Proof:} i) Let $(I,J,K)$ be the hypercomplex trivialization of $S^2 H$ on $(M,g)$.
It is clear that $I$, $J$, $K$ are orthogonal with respect to $g_0$. The
K\"ahler form of $I$ with respect to $g_0$ is
$$\Omega_I^0 = \frac{1}{\nu ^2} ( g(\varphi,\varphi ) + \nu )\left( \varphi \wedge
I\varphi + J\varphi \wedge K\varphi + \nu \Omega_I \right).$$
Then, using \pref{60} and \pref{58}, it is easily verified that
$\frac{1}{2}dId\mu = \Omega_I^0$, and similarly for $J$ and
$K$. Hence, if $g_0$ is non-degenerate, then, by
Lemma~\ref{lem1}~ii), it is pseudo-hyper-K\"ahler with
hyper-K\"ahler potential $\mu$. On the orthogonal complement
of $span_\mathbb{H} \{ \xi \}$ we have $g_0 =
\frac{1}{\nu}(g(\varphi ,\varphi) + \nu) g$ and $g_0 (\xi,
\xi) = \frac{1}{\nu^2}(g(\varphi ,\varphi) + \nu)^2
g(\varphi ,\varphi)$, which is positive if $\varphi \not =
0$. This proves the non-degeneracy of $g_0$ and the
assertion about its signature.

ii) Obviously, the given hypercomplex structure and $g_p$ form an
almost quaternionic Hermitian structure. The K\"ahler form of $I$ with respect
to $g_p$ is
$$\Omega_I^p =  - \frac{p}{(pg_0(d\mu,d\mu ) + 1)^2 } (d\mu \wedge Id\mu + Jd\mu
\wedge Kd\mu ) + \frac{1}{pg_0(d\mu,d\mu ) + 1} \Omega_I^0$$
and similarly for $J$ and $K$. Now, using \pref{67} and \pref{65}, it is easy to
see that $\Omega_I^p$, $\Omega_J^p$, $\Omega_K^p$ satisfy \pref{21} for $a =
I\varphi_p$, $b = J\varphi_p$,  $c = K\varphi_p$. Thus, by
Lemma~\ref{lem1}~i) and Proposition~\ref{prop1}, $(M,g_p)$ is a hHqK
manifold with Lee form $\varphi_p$. The statement about the reduced scalar
curvature follows from Lemma~\ref{2} by a straightforward computation.

Let $\zeta$ be the vector field dual to $d\mu$ with respect to $g_0$. Then $g_p
(\zeta, \zeta) = \frac{g_0 (d\mu ,d\mu)}{(pg_0(d\mu,d\mu ) + 1)^2}$ and hence
if $g_0$ is positive definite on $span_\mathbb{H} \{ d\mu \}$, then so is
$g_p$. On the orthogonal complement of $span_\mathbb{H} \{ \zeta\}$ we have
$g_p = \frac{1}{pg_0(d\mu,d\mu ) + 1}g_0$. This completes the proof of the
non-degeneracy of $g_p$ and the statement about its positive definiteness.

Part iii) is straightforward, after noticing that $g(\varphi
,\varphi) = \frac{1}{4} \nu^2 g_0(d\mu,d\mu )$.  \hfill $\square$

\vspace{3mm}
We illustrate the above theorem by the following simplest example.

\vspace{3mm}
\noindent
{\bf Example 1} Let
$g_0 = {\rm Re}
\left( \sum _{\lambda = 1}^n d\bar{x}^\lambda
\otimes dx^\lambda \right)$
be the standard flat metric on $\mathbb{H}^n$. Then
$(\mathbb{H}^n ,g_0 )$ is a hyper-K\"ahler manifold with hyper-K\"ahler
potential $\mu = \frac{1}{2}\Vert x\Vert ^2$, the hypercomplex structure being
the standard one of $\mathbb{H}^n$ (defined by \pref{2} after the standard
identification of the tangent spaces of $\mathbb{H}^n$ with $\mathbb{H}^n$).

Now let us consider the basic examples of quaternionic K\"ahler manifolds:
the quaternionic projective space $\mathbb{H}P^n = Sp(n+1)/Sp(n)\times Sp(1)$
and its dual symmetric space $\mathbb{H}H^n = Sp(1,n)/Sp(n)\times Sp(1)$.

An equivalent definition of $\mathbb{H}P^n$ is $\mathbb{H}P^n =
(\mathbb{H}^{n+1}\backslash \{0\})/\sim$, where for $x, y \in
\mathbb{H}^{n+1}\backslash \{0\}$ we have $x \sim y$ iff $y = xq$ for some $q
\in \mathbb{H}\backslash \{0\}$. Let $U_0 = \{[x^0,\dots,x^n] \in \mathbb{H}P^n
: x^0 \not = 0\}$. Then $U_0$ is a domain of non-homogeneous quaternionic
coordinates on $\mathbb{H}P^n$, which are given by
$$U_0 \ni
[1,x^1,\dots,x^n] \longleftrightarrow x=(x^1,\dots,x^n)\in \mathbb{H}^n.$$
Thus $U_0$ is naturally diffeomorphic to $\mathbb{H}^n$ and, moreover, the
quaternionic structure on $U_0$ is spanned by the standard hypercomplex
structure of $\mathbb{H}^n$.

Similarly, $\mathbb{H}H^n$ is diffeomorphic to $\{x \in \mathbb{H}^n : \Vert x\Vert < 1\}$ and again the quaternionic structure comes from the standard hypercomplex
structure of $\mathbb{H}^n$.

Hence, $U_0 \subset \mathbb{H}P^n$ and $\mathbb{H}H^n$ are hHqK manifolds. In
the given coordinates their metrics are
$$g_{\pm} = \frac{1}{(1 \pm \Vert x\Vert ^2
)^2} {\rm Re} \left(\left(1 \pm \Vert x\Vert ^2 \right) \sum _{\lambda = 1}^n
d\bar{x}^\lambda \otimes dx^\lambda \mp \left(\sum _{\lambda = 1}^n
d\bar{x}^\lambda . x^\lambda \right) \otimes \left(\sum _{\lambda = 1}^n
\bar{x}^\lambda . dx^\lambda \right)\right).$$

The hHqK manifolds  $(U_0 ,g_+)$ and $(\mathbb{H}H^n ,g_-)$ are obtained
from $(\mathbb{H}^n ,g_0 )$ by the construction in part ii) of Theorem~\ref{th2}
for parameters $p=1$ and $p=-1$ respectively. In particular, their Lee forms are
given by $\varphi_\pm = - d \ln (1 \pm \Vert x\Vert ^2 )$ (and the reduced
scalar curvatures are $\nu_\pm = \pm 4$).

\vspace{3mm}
Now we summarize some results of Swann \cite{Sw}.

Let $M'$ be a quaternionic K\"ahler manifold and $P'$ be the principal
$SO(3)$-bundle over $M'$, whose points are the frames $(I',J',K')$ trivializing
$S^2 H$ and satisfying \pref{3}. The {\it Swann bundle} over $M'$ is the
principal $\mathbb{R}_+ \times SO(3)$-bundle $\mathcal{U}(M') = \mathbb{R}_+
\times P'$. The Levi-Civita connection defines a horizontal distribution on
$P'$ and hence also on $\mathcal{U}(M')$. A hypercomplex structure $(I,J,K)$ is
defined on $\mathcal{U}(M')$ in the following way. The projection $\pi :
\mathcal{U}(M') \longrightarrow M'$ induces an isomorphism of the horizontal
space on $\mathcal{U}(M')$ at the point $(r,I',J',K')$ and the tangent space
of $M'$ at the corresponding point. On the horizontal space  $I$, $J$, $K$ are
defined to correspond respectively to $I'$, $J'$, $K'$ under this isomorphism.
On the fibres $(I,J,K)$ is the standard hypercomplex structure, that is, $I$,
$J$, $K$ are given by \pref{2} after identifying the tangent spaces of
$\mathbb{R}_+ \times SO(3) = \mathbb{H}^*/{\mathbb{Z}_2}$ with its Lie algebra
$\mathbb{H}$.

There exist also quaternionic K\"ahler metrics compatible with this
hypercomplex structure.
\begin{theo}\label{th3}
\cite{Sw} Let $(M',g')$ be a $4(n-1)$-dimensional quaternionic K\"ahler manifold
with reduced scalar curvature $\nu '$ and $r$ be the radial coordinate on
$\mathcal{U}(M')$. Then for $p \not =0$ the above
hypercomplex structure and \begin{equation}\label{76}
g_p = \frac{1}{(pr^2 + 1)^2 } (\pmb{1} + \pmb{L}) \, dr \otimes dr +
\frac{\nu ' r^2}{4(pr^2 + 1)} \pi^* g'
\end{equation}
form on the
submanifold $\nu ' (pr^2 +1) >0$ a (positive definite) hHqK structure with Lee
form $\varphi _p = - d\ln \vert pr^2+1 \vert$ and reduced scalar curvature $\nu
_p = 4p$. The metric $g_0$ is pseudo-hyper-K\"ahler with respect to the same
hypercomplex structure and has hyper-K\"ahler potential $\mu = \frac{r^2}{2}$.
It has Riemannian signature if $\nu ' >0$ and signature $(4,4(n-1))$ if $\nu '
<0$. \end{theo}
Notice that $(\mathcal{U}(M'),g_p)$ and  $(\mathcal{U}(M'),g_q)$ are homothetic if $p$ and $q$ have
the same sign.

It is easily seen  that $g_p (\varphi_p ,\varphi_p) + \nu_p = 4p(pr^2 +1)$.
Thus, Theorem~\ref{th3} gives examples of hHqK manifolds with exact Lee form
$\varphi$ such that $\varphi$ is everywhere non-zero and $\vert \varphi \vert
^2 +\nu \not = 0$. It turns out, as suggested by Theorem~\ref{th2}, that these
examples exhaust locally all such manifolds.
\begin{theo}\label{th4}
Let $(M,g)$ be a hHqK manifold with closed Lee form $\varphi$ and reduced
scalar curvature $\nu$ such that $\varphi \not = 0$ and $\vert \varphi \vert ^2
+\nu \not = 0$. Then $(M,g)$ is locally isometric to $(\mathcal{U}(M'),
g_{\frac{1}{4}\nu})$ for some quaternionic K\"ahler manifold $M'$.
\end{theo}

\noindent {\it Proof:} From the results of Swann \cite{Sw} it follows that
every pseudo-hyper-K\"ahler manifold with hyper-K\"ahler potential $\mu$, such
that $d\mu$ does not vanish, is locally homothetic to $(\mathcal{U}(M'),g_0)$ for some
pseudo-quaternionic K\"ahler manifold $M'$. From \pref{76} we get
$$g_p = -\frac{pr^2}{(pr^2 + 1)^2 } (\pmb{1} + \pmb{L}) \, dr \otimes dr +
\frac{1}{pr^2 + 1} g_0.$$
Since the hyper-K\"ahler potential of
$g_o$ is $\mu = \frac{r^2}{2}$, we obtain the result by applying
Theorem~\ref{th2}.   \hfill $\square$

\vspace{3mm}
Next, we consider the case when $\varphi$ vanishes at some point.

The condition $R' \equiv 0$ ensures that the Cauchy problem for \pref{56}
locally has a solution for any  initial data. Therefore, on $\mathbb{H}P^n$ and
$\mathbb{H}H^n$ the bundle $S^2 H$ can be locally trivialized by a
hypercomplex structure, whose Lee form is closed and vanishes at some point. In
fact, these are the only such manifolds:
\begin{theo}\label{th5}
i) A hHqK manifold, whose Lee form is closed and vanishes at some point, is
locally homothetic to $\mathbb{H}P^n$ or $\mathbb{H}H^n$.

ii) A hyper-K\"ahler manifold with hyper-K\"ahler potential $\mu$, such that
$d\mu$ vanishes at some point, is flat.
\end{theo}

\noindent {\it Proof:} i) Differentiating \pref{57} and
using also \pref{56}, we see that
\begin{equation}\label{78} 
\nabla_X R'(\xi,Y,Z,W) = \frac{1}{2} \nu R'(X,Y,Z,W).
\end{equation}
Thus, at the point $p$, where $\varphi$ vanishes, we have $R' =0$.

We shall prove that all the covariant derivatives of the curvature
tensor also vanish at $p$. Then, because the metric is Einstein and hence
analytic, it will follow that the curvature is $R = \frac{1}{4} \nu R_0$, that
is, the manifold is locally homothetic to $\mathbb{H}P^n$ or $\mathbb{H}H^n$.

Using
\pref{56}, \pref{57} and \pref{78}, it is easily proved by induction that
\begin{eqnarray}\label{79}
\nabla^{k+1}_{X_1,\dots,X_{k+1}} R'(\xi,Y,Z,W) = \frac{1}{2} \nu
\sum_{s=1}^{k+1} \nabla^k_{X_1,\dots,\widehat{X}_s,\dots,X_{k+1}}
R'(X_s,Y,Z,W) && \\
+ P_k(\varphi,R',\dots,\nabla^{k-1}R')(X_1,\dots
,X_{k+1},Y,Z,W), && \nonumber
\end{eqnarray}
where $P_k$ is a polynomial such that $P_k(\cdot,0,\dots,0) = 0$ (in
particular, $P_k$ has no term of order zero). The notation
$\widehat{X}_s$ is used to indicate that the argument $X_s$ is omitted. Now,
supposing that $\nabla^l R' =0$ at $p$ for $l<k$, we see by \pref{79} that
\begin{equation}\label{80} \sum_{s=1}^{k+1}
\nabla^k_{X_1,\dots,\widehat{X}_s,\dots,X_{k+1}} R'(X_s,Y,Z,W)(p) = 0.
\end{equation}
Since the antisymmetrization of $\nabla^k_{X_1,\dots,X_k}  R'(X_{k+1},Y,Z,W)$
with respect to $X_1$ and $X_2$ is $(R(X_1,X_2)\nabla^{k-2}
R')(X_3,\dots,X_{k+1},Y,Z,W)$, it follows that at $p$ it is symmetric with
respect to $X_1$ and $X_2$. Similarly, it is symmetric at $p$ with respect to
$X_s$ and $X_{s+1}$ for every $s<k$, since its antisymmetrization with respect
to these two arguments is  expressed by the covariant derivatives of $R'$ of
order less than $k-1$.

Hence, $\nabla^k R' (p)$ is symmetric with respect to the first $k$ arguments
and the proof of i) is completed by the following algebraic lemma:

\begin{lem}\label{lem3}
Let $T \in S^k T^* \otimes \Lambda^2 T^*$ satisfy
the Bianchi identity with respect to the last three arguments and
\begin{equation}\label{82}
\sum_{s=1}^{k+1} T(X_1,\dots,\widehat{X}_s,\dots,X_{k+1},X_s,X_{k+2}) = 0.
\end{equation}
Then $T=0$.
\end{lem}

\noindent {\it Proof:} Antisymmetrizing \pref{82} with respect to $X_{k+1}$
and $X_{k+2}$ and using the Bianchi identity with respect to
the last three arguments, we obtain
$$2T(X_1,\dots,X_{k+2}) + \sum_{s=1}^k
T(X_1,\dots,\widehat{X}_s,\dots,X_k,X_s,X_{k+1},X_{k+2}) = 0.$$
The symmetry with respect to the first $k$ arguments now implies
$T=0$. \hfill $\square$

\vspace{3mm}
ii) It follows from \pref{66} that $R(X,Y,Z,\zeta) = 0$, where $\zeta =
(d\mu)^\#$. Now ii) can be proved in the same way as i), the polynomials $P_k$
being identically zero.  \hfill $\square$

\vspace{3mm}
\noindent {\bf Remarks }

{\bf 6)} Lemma~\ref{lem3} in fact verifies that the condition \pref{80} forces
the vanishing of the component of $\nabla^k R'$ in the subspace of highest
dominant weight in the space of tensors with the symmetries of the $k$th
covariant derivative of the curvature tensor of an Einstein manifold. Thus
Theorem~\ref{th5} follows from Theorem~10.2 in \cite{S2} or from the similar
result for the special case of a quaternionic K\"ahler manifold, given in the
proof of Theorem~2.6 in \cite{PS}.

{\bf 7)} Part ii) of Theorem~\ref{th5} also follows from the results in
\cite{A2,Y}. By \pref{66}, we see that $\zeta$ is an infinitesimal conformal
transformation with non-vanishing divergency. If $d\mu$ vanishes somewhere,
then $\zeta$ is an essential infinitesimal conformal transformation, that is,
it is not an infinitesimal isometry for any conformal metric. Then it follows
from the "obvious" parts of Proposition~2 in \cite{A2} or the Theorem in
\cite{Y} that the manifold is locally conformally flat (notice that in this
obvious part it is not necessary to have a global conformal transformation).
But it is also Ricci-flat, and hence flat.

{\bf 8)} Part i) of Theorem~\ref{th5} can also be proved using part ii) and
Theorem~\ref{th2}: Since $\nu (\vert \varphi \vert ^2 + \nu ) >0$ in a
neighbourhood of the point in which $\varphi =0$, the
hyper-K\"ahler metric given by Theorem~\ref{th2}~i) is positive definite. Also,
its hyper-K\"ahler potential $\mu$ attains its minimum at this point and hence
$d\mu =0$ at it. Now part i) of Theorem~\ref{th5} follows from part ii) and
Example~1.

\vspace{3mm}
Finally, we focus our attention on the case of hHqK manifold with closed Lee
form $\varphi$, satisfying $\vert \varphi \vert ^2 + \nu = 0$. Clearly, this is
possible only if $\nu <0$. First we construct examples of such manifolds.

Let $(M',g')$ be a $4(n-1)$-dimensional hyper-K\"ahler manifold, the
hypercomplex structure being $(I',J',K')$. Then the K\"ahler
forms $\Omega_{I'}$, $\Omega_{J'}$, $\Omega_{K'}$ are closed. Thus, restricting
our considerations on a small open set, we can fix 1-forms $\alpha_{I'}$,
$\alpha_{J'}$, $\alpha_{K'}$ such that
$$\Omega_{I'} = d\alpha_{I'}, \quad  \Omega_{J'} = d\alpha_{J'}, \quad
\Omega_{K'} = d\alpha_{K'}.$$
Let $(t,u,v,w)$ be the standard coordinates on $\mathbb{H} \cong
\mathbb{R}^4$, that is,  $\mathbb{H} \ni q = t + ui + vj + wk$. Let $M = \{
q \in \mathbb{H}\, :\, t>0 \} \times M'$ and $\pi: M \longrightarrow M'$ be
the projection. We fix a negative constant $\nu$ and define on $M$ an almost
hypercomplex structure $(I,J,K)$ by
\begin{equation}\label{88} 
Idt = -du + \nu \pi^*\alpha_{I'}, \quad  Jdt = -dv + \nu \pi^*\alpha_{J'},
\quad  Kdt = -dw + \nu \pi^*\alpha_{K'},
\end{equation}
\begin{equation}\label{89}
I\pi^* \beta = \pi^* I'\beta, \quad  J\pi^* \beta = \pi^* J'\beta, \quad
K\pi^* \beta = \pi^* K'\beta, \qquad \beta \in T^* M'
\end{equation}
and a Riemannian metric $g$ by
\begin{equation}\label{90} 
g = -\frac{1}{\nu t^2}(\pmb{1} + \pmb{L}) \, dt \otimes dt + \frac{1}{t}
\pi^* g'.
\end{equation}
\begin{theo}\label{th6}
The manifold $(M,g,I,J,K)$ is hHqK with reduced scalar curvature $\nu$
and Lee form $\varphi = -d\ln t$, which satisfies $\vert \varphi
\vert ^2 + \nu = 0$.
\end{theo}

\noindent{\it Proof:} It is clear from the definitions that $I$, $J$, $K$ are
orthogonal with respect to $g$. The K\"ahler form of $I$ is
$$\Omega_I = -\frac{1}{\nu t^2} \bigl( dt \wedge (-du + \nu \pi^*\alpha_{I'}) +
(-dv + \nu \pi^*\alpha_{J'}) \wedge  (-dw + \nu \pi^*\alpha_{K'}) \bigr) + \frac{1}{t}
\pi^* \Omega_{I'}$$
and similarly for $\Omega_J$ and $\Omega_K$. Now a straightforward computation
shows that \pref{21} is satisfied with
$a = -Id\ln t$, $b = -Jd\ln t$, $c = -Kd\ln t$.
Thus, by Lemma~\ref{lem1} and Proposition~\ref{prop1}, $(g,I,J,K)$ is a hHqK
structure on $M$ with Lee form $\varphi = -d\ln t$. That the
reduced scalar curvature is $\nu$ is verified using Lemma~\ref{lem2}. The
equality $\vert \varphi \vert ^2 + \nu = 0$ is an obvious consequence of the
definition of $g$.   \hfill $\square$

\vspace{3mm}
\noindent {\bf Remark 9} It is easy to see that
$I\frac{\partial}{\partial t} = - \frac{\partial}{\partial u}$,
$J\frac{\partial}{\partial t} = - \frac{\partial}{\partial v}$,
$K\frac{\partial}{\partial t} = - \frac{\partial}{\partial w}$,
that is, on the fibres of $\pi$ the hypercomplex structure is the standard
hypercomplex structure on $\mathbb{H}$. Also, $\pi : (M,tg)
\longrightarrow (M',g')$ is a Riemannian submersion.

\vspace{3mm}
Now we show that locally the converse of Theorem~\ref{th6} is also true.
\begin{theo}\label{th7}
A hHqK manifold $M$ with closed Lee form $\varphi$ and negative reduced scalar
curvature $\nu$, such that $\vert \varphi \vert ^2 + \nu = 0$, is locally
isometric to one of the manifolds in Theorem~\ref{th6}.
\end{theo}

\noindent{\it Proof:} Since $\nabla \psi$ is symmetric, the distribution $\eta
^\bot = \{ X  : g(\eta, X) = 0 \}$ is integrable. Furthermore, the integral
curves of $\eta$ are geodesics (up to a change of the parameter) and therefore
$X \in \eta^\bot$ implies $[\eta,X] \in \eta^\bot$. Hence, in a neighbourhood
of a fixed point $p_0 \in M$ we can choose coordinates
$(t,u,v,w,x^1,\dots,x^{4(n-1)})$ such that
\begin{equation}\label{94} 
\eta = -\frac{\partial}{\partial t}, \quad \eta^\bot = span \left\{
\frac{\partial}{\partial u},\frac{\partial}{\partial
v},\frac{\partial}{\partial w},\frac{\partial}{\partial
x^1},\dots,\frac{\partial}{\partial x^{4(n-1)}} \right\}.
\end{equation}
By \pref{64}, we have that $\eta$, $I\eta$, $J\eta$, $K\eta$ commute ($C=0$ by
\pref{61}) and therfore we can take the above coordinates in such a way that
\begin{equation}\label{95} 
I\eta = \frac{\partial}{\partial u}, \quad
J\eta = \frac{\partial}{\partial v}, \quad
K\eta = \frac{\partial}{\partial w}.
\end{equation}

From now on we restrict our considerations to this coordinate neighbourhood.

Since  $\vert \varphi \vert ^2 + \nu =0$, we have  $\vert \psi \vert ^2 = -\nu
e^{2f}$, and it follows from \pref{94} that $\psi = \nu e^{2f} dt$. On the
other hand, $\psi = d e^f$ and therefore $e^f = \frac{1}{-\nu t +D}$, where $D$
is a constant. Changing the coordinate $t$ by a translation, we can assume that
$e^f = -\frac{1}{\nu t}$ (and hence $t>0$).

Let $M' = \{ p \in M : t(p) = t(p_0),\, u(p) = u(p_0),\, v(p) = v(p_0),\, w(p) =
w(p_0) \}$ and $\pi : M \longrightarrow M'$ be the projection. We call
$\mathcal{V} = span_\mathbb{H} \{ \eta \}$ {\it the vertical distribution} and
its orthogonal complement $\mathcal{H}$ {\it the horizontal distribution}. $\mathcal{V}$
and $\mathcal{H}$ are invariant under the action of the complex structures  $I$, $J$,
$K$. As seen before, since $\vert \varphi \vert ^2 + \nu = 0$, each of $\eta$,
$I\eta$, $J\eta$, $K\eta$ is an infinitesimal hypercomplex automorphism.
Therefore $I$, $J$, $K$ project down to almost complex structures $I'$, $J'$,
$K'$ on $M'$, that is, \pref{89} is satisfied.

Let $h: TM \longrightarrow \mathcal{H}$ be the orthogonal projection. We define a
Riemannian metric $g'$ on $M'$ by
$$g'(X,Y) = tg(hX,hY), \qquad X,Y \in T_p M' \subset T_p M.$$
It is not difficult to see that the tensor $tg(h\cdot,h\cdot)$ on $M$
projects down to a tensor on $M'$, that is, $\pi : (M,tg) \longrightarrow
(M',g')$ is a Riemannian submersion. Thus, $g$ is given by \pref{90}.

It is obvious that $I'$, $J'$, $K'$ are orthogonal with respect to $g'$. So,
\begin{equation}\label{98} 
\Omega_I = -\frac{1}{\nu t^2} \bigl( dt \wedge Idt + Jdt \wedge  Kdt \bigr) +
\frac{1}{t} \pi^* \Omega_{I'}.
\end{equation}
Now, it follows from \pref{95} that
\begin{equation}\label{99} 
Idt = -du + \sum_{s=1}^{4(n-1)} f_s dx^s.
\end{equation}
Hence,
\begin{equation}\label{100} 
dIdt = \sum_{s=1}^{4(n-1)} df_s \wedge dx^s.
\end{equation}
But $\varphi = df$ and therefore $dt = - t\varphi$. Thus, using
\pref{58} and \pref{98}, we obtain
\begin{equation}\label{101}
dIdt = \nu \pi^* \Omega_{I'}.
\end{equation}
It follows by \pref{100} and \pref{101} that the coefficients $f_s$ do not
depend on $t$, $u$, $v$, $w$ and therefore
\begin{equation}\label{102}
\sum_{s=1}^{4(n-1)} f_s dx^s = \nu \pi^* \alpha_{I'}.
\end{equation}
for some form  $\alpha_{I'}$ on $M'$. Now, by \pref{99}, \pref{101} and
\pref{102}, we obtain $\pi^* d\alpha_{I'} = \pi^* \Omega_{I'}$, that is,
$d\alpha_{I'} = \Omega_{I'}$. Repeating the same argument for $J$ and $K$, we
see by Lemma~\ref{lem1} ii) that $M'$ is hyper-K\"ahler and \pref{88} is
satisfied.   \hfill $\square$

\vspace{3mm}
\noindent{\bf Remark 10} The above proof can be easily modified when $\vert
\varphi \vert ^2 + \nu \not = 0$ to get a proof of Theorem~\ref{th3}. In this
case we can again take $\eta = -\frac{\partial}{\partial t}$. The vector fields
$I\eta$, $J\eta$, $K\eta$ do not commute, but they form an integrable
distribution and the coordinates can be taken so that $span \{ I\eta, J\eta,
K\eta \} = span \{ \frac{\partial}{\partial u},\frac{\partial}{\partial
v},\frac{\partial}{\partial w} \}$. Again we have a Riemannian submersion  $\pi
: (M,\frac{1}{e^f \vert \varphi \vert ^2}g) \longrightarrow (M',g')$, where
$$g' = \frac{1}{e^f \vert \varphi \vert^2} g(h\cdot,h\cdot)_{\vert_{M'}} =
\frac{1}{C - \nu e^f} g(h\cdot,h\cdot)_{\vert_{M'}}.$$
The complex structures $I$, $J$, $K$ do not project down to $M'$, but their
span does and together with $g'$ forms a quaternionic K\"ahler structure with
reduced scalar curvature $C\nu$. In the chosen coordinates
$$e^f = \frac{C}{e^{Ct} + \nu }.$$
Now, after changing the coordinate $t$ by $e^{Ct} = \frac{1}{4}\nu^2 r^2$, it is
easily seen that $g$ is locally isometric to the metric $g_{\frac{1}{4}\nu}$ on
$\mathcal{U}(M')$ in Theorem~\ref{th3}.

\section{Complete hHqK manifolds}\label{sec5}

It is proved in \cite{AM2,AMP} that on a compact quaternionic K\"ahler manifold
the bundle $S^2 H$ cannot be globally trivialized by an almost hypercomplex
structure. In particular, this is true for the complete quaternionic K\"ahler
manifolds with positive scalar curvature.
It is conjectured by Alekseevsky, Marchiafava and Pontecorvo \cite{AMP} that
the only complete simply connected hHqK manifold is $\mathbb{H}H^n$. In support
of this they prove that if in addition the Lee form is closed, then there
exists a (possibly singular) integrable quaternionic distribution, whose
regular orbits are locally homothetic to $\mathbb{H}H^k$. Below we show that
under this additional assumption the manifold is indeed $\mathbb{H}H^n$.
\begin{prop}\label{prop4}
A complete simply connected hHqK manifold with
closed Lee form is homothetic to $\mathbb{H}H^n$.
\end{prop}

\noindent {\it Proof:} If $\varphi$ vanishes somewhere, then by
Theorem~\ref{th5} the manifold is homothetic to $\mathbb{H}H^n$ (as remarked
above, $\nu < 0$). Thus, it is enough to consider the case when $\varphi \not =
0$ everywhere.

Let $x_t$ be a geodesic parametrized with respect to its length $t$ and
tangent at some point to $\xi$. Then, because of the completeness $x_t$ is
defined for all $t \in \mathbb{R}$ and by \pref{106} it is tangent to $\xi$ at
each point. Thus, $\xi_{x_t} = h(t) {\dot{x}}_t$ with $h(t) = \vert \xi_{x_t}
\vert \not = 0$ and \pref{106} becomes
$$\frac{dh}{dt} = - \frac{1}{2} (h^2 + \nu ).$$
This equation has no solutions
with $h^2 + \nu > 0$, defined on the whole $\mathbb{R}$, while every solution
with $h^2 + \nu < 0$ vanishes somewhere.

Hence, it remains to consider the case $\vert \varphi \vert^2 + \nu =0$. As
seen before, the distribution $span_{\mathbb{H}} \{ \eta \}$ is totally
geodesic and its integral manifolds are of constant (negative) curvature $\nu$.
Every quaternionic K\"ahler manifold is analytic and hence, by Proposition~7 in
\cite{AM3}, these totally geodesic submanifolds can be extended to complete
(immersed) totally geodesic submanifolds. By \pref{64}, we have that $I\eta$,
$J\eta$, $K\eta$ are three commuting Killing vector fields on them. This is a
contradiction, since the maximal commutative subalgebra of the algebra
$\mathfrak{so}(1,4)$ of Killing vector fields of $\mathbb{R}H^4$ is
2-dimensional. \hfill $\square$

\vspace{3mm}
\noindent{\bf Remark 11}
Using \pref{66}, it can be easily proved in the same way as above that on a
complete simply connected hyper-K\"ahler manifold with hyper-K\"ahler potential
$\mu$ there exists a point at which $d\mu$ vanishes. Hence, by
Theorem~\ref{th5}~ii), the only such manifold is $\mathbb{H}^n$ with the flat
metric.

\vspace{3mm} The proof of Proposition~\ref{prop4} shows that there are
no complete hHqK manifolds whose Lee form satisfies $\vert \varphi \vert^2 +
\nu =0$. On the other hand, there exists a global solution of \pref{56} on
$\mathbb{H}H^n$, satisfying this condition. This follows from the local
existence on $\mathbb{H}H^n$ of a solution of the Cauchy problem for \pref{56}
with any initial data. If $\vert \varphi \vert^2 + \nu =0$ at one point, then
this is true everywhere, where the solution is defined, and since the isometry
group of $\mathbb{H}H^n$ acts transitively on the unit tangent bundle, this
local solution can be extended on the whole $\mathbb{H}H^n$. The corresponding
vector field $\eta$ is an infinitesimal quaternionic automorphism which is not
a Killing vector field. It can not be a complete vector field since the group
of quaternionic automorphisms of $\mathbb{H}H^n$ coincides with the group of
its isometries \cite{AM3}.

It follows from the above discussion that on  $\mathbb{H}H^n$ the bundle
$S^2 H$ can be locally trivialized by hypercomplex structures which have the
same (globally defined) Lee form $\varphi$. Such a situation cannot occur on a
compact quaternionic K\"ahler manifold:

\begin{prop}\label{prop5} On a compact quaternionic K\"ahler manifold the
equation \pref{25} has no global solutions.
\end{prop}

\noindent {\it Proof:} By Remark~4, the exact form $\Phi$ is
harmonic and therefore $\Phi = 0$. Now \pref{60} shows
that at every critical point of $\vert \varphi \vert^2$ we have $(\vert \varphi
\vert^2 + \nu) \varphi =0$. Integrating \pref{44}, we see that $\nu >0$.
Hence, at a point of maximum of $\vert \varphi \vert^2$ the form $\varphi$ must
vanish. Thus, $\varphi \equiv 0$ and therefore $\nu = 0$, which is a
contradiction. \hfill $\square$

\vspace{3mm}
The quaternionic projective space is the only complete locally hHqK
manifold with positive scalar curvature in dimensions 4 and 8. This follows
from Theorem~\ref{th1} and the results of
Hitchin \cite{H2}, Friedrich and Kurke \cite{FK} and Poon and Salamon
\cite{PS}, according to which every complete quaternionic K\"ahler manifold
with positive scalar curvature in these dimensions is symmetric. In fact, there
are no known examples of non-symmetric complete quaternionic K\"ahler manifold
with positive scalar curvature. Thus, it seems reasonable to expect that
$\mathbb{H}P^n$ could be characterized by the above property in all dimensions.

\vspace{3mm}
\noindent {\it Acknowledgements.} I would like to thank Paul Gauduchon for
drawing my attention to the problem discussed in the present paper during a
lecture given by him in Berlin in 1999. I am also grateful to Gueo
Grantcharov, Vestislav Apostolov, Florin Belgun and Volker Buchholz for some
useful remarks and discussions.

\vspace{10mm}
\noindent
Bogdan Alexandrov \\
Humboldt Univesit\"at zu Berlin \\
Institut f\"ur Mathemathik \\
Sitz: Rudower Chaussee 25  \\
10099 Berlin \\
{\tt e-mail: \quad bogdan@mathematik.hu-berlin.de}

\end{document}